\documentclass[12pt,twoside,sort&compress ]{article} 
\usepackage{mathrsfs,amsfonts,amsmath,amssymb,amsbsy,tikz}\usepackage{fontenc}\usepackage{textcomp}
\usepackage{graphicx}
\usepackage{amssymb}
\usepackage{graphicx}

\def\ge{\geq}

\def\bbb{\begin{eqnarray*}}

\def\eee{\end{eqnarray*}}

\begin{document} 
\baselineskip=18pt
\begin{center}
\vspace{-0.6in}{\large \bf Devaney chaos in non-autonomous
discrete systems }
\\ [0.2in]
 Hao Zhu,\ \ Yuming Shi$^{\dag}$, \ \ Hua Shao\\
\vspace{0.15in}  Department of Mathematics, Shandong University\\
Jinan, Shandong 250100, P. R. China\\

\footnote{$^\dag$ The corresponding author}
\footnote{Email addresses: haozhu@mail.sdu.edu.cn(H. Zhu), ymshi@sdu.edu.cn(Y. Shi), huashaosdu@163.com (H. Shao).}
\footnote{This research was supported by the NNSF of China (Grant 11571202).}
\end{center}

{\bf Abstract.}  This paper is concerned with   Devaney chaos in non-autonomous discrete  systems.
  It is shown that in its definition, the two former conditions, i.e., transitivity and density of periodic points,
in a set imply the last one, i.e., sensitivity, in the case that the set is unbounded, while a similar result holds
under two additional conditions in the other case that the set is bounded. Furthermore,
some chaotic behavior is studied for a class of non-autonomous systems, each of which is governed by a convergent sequence of continuous maps.
\medskip

\medskip

{2010 {\bf \it Mathematics Subject Classification}}: 37B55; 37D45;
54H20.

\bigskip

\noindent{\bf 1. Introduction}\medskip

In this paper, we consider the following non-autonomous discrete  system:
\vspace{-0.2cm}$$x_{n+1}=f_n(x_n),\;\;\;\; n\geq 0,                             \eqno(1.1)\vspace{-0.2cm}$$
where $f_n: D_n\to D_{n+1}$ is a map and $D_n$
is a subset of a metric space $(X, d)$ for each $n\geq0$.
When $f_n=f$ and $D_n=D$ for each $n\geq0$, system (1.1) becomes the following
autonomous discrete  system:
\vspace{-0.2cm}$$x_{n+1}=f(x_n),\;\;\;\; n\geq 0.                               \eqno(1.2)\vspace{-0.2cm}$$

Chaos, which is one of the central topics of research on nonlinear science,
reflects the complexity of dynamical systems.
Chaos theory of autonomous systems has been extensively studied [7, 17, 25] and
there appears several definitions   of chaos in autonomous systems
[5, 7, 15, 18, 25]. In particular, Devaney in 1989 gave the following definition of chaos
for  system (1.2):\medskip

\noindent{\bf Definition 1.1} [7]. Let $X$ be a metric space. System (1.2) is said to be chaotic on $X$ if \begin{itemize}\vspace{-0.2cm}
\item[{\rm (1)}] system (1.2) is topologically transitive in $X$;\vspace{-0.2cm}
\item[{\rm (2)}] the periodic points of system (1.2) is dense in $X$;\vspace{-0.2cm}
\item[{\rm (3)}] system (1.2) has sensitive dependence on initial conditions in $X$.\vspace{-0.15cm}
\end{itemize}
In 1992, Banks and his coauthors proved that (1) and (2)  imply (3) in the above definition
 if $f$ is continuous in $X$ [3]. This shows that condition (3) in this definition is
redundant in this case. For further results, see [10, 12].

Non-autonomous systems occur more often than autonomous ones since many physical,
biological, and economical problems  are often  described  by non-autonomous systems.
Thus, the study on complexity of system (1.1) has been of  increasing interest by many
researchers in recent years [2, 4, 6, 8, 11, 14, 16, 19--24, 26]. We briefly recall
some existing results related to Devaney chaos for system (1.1). In 2006, Tian and Chen
studied system (1.1) in the special case that $D_n=X$ for each $n\geq 0$ and extended
the concept of  Devaney chaos to system (1.1) [24].  In 2009, Shi and Chen
removed the assumption that $D_n=X$ for each $n\geq 0$ and thus generalized
the concept of Devaney chaos for  system (1.1) [21] (See also Definition 2.4).
There appears a natural question:  whether the  two former conditions--transitivity
and density of periodic points  imply the last one--sensitivity for  system (1.1).
This is  the main topic studied  in the present  paper.

Chaotic behavior of system (1.1) was  studied
 in the special case  that $\{f_n\}_{n=0}^{\infty}$ converges to a map $f$ by some scholars (cf., [6, 8, 26] and references cited therein).  C\'{a}novas investigated the relationships between Li-Yorke chaos of system (1.1) and that of system (1.2) for interval maps [6]. Dvo${\rm\check{r}\acute{a}}$kov${\rm\acute{a}}$ showed that system (1.2) may not be distributionally chaotic even if system (1.1) has a distributionally scrambled set with
full Lebesgue measure [8]. In the present paper, we shall study Devaney chaos of system (1.1)  in this special case.

The rest of this paper is organized as follows.
In Section 2, some basic concepts and lemmas are given. In particular,
properties of the three conditions of Devaney chaos for  system (1.1) are investigated,
some of which are shown by different methods from those in the autonomous case.
In Section 3,   it is shown that
the two former  conditions in a set $V$ imply the last one of Devaney chaos
for system (1.1) in the case that $V$ is unbounded,
while  a similar result holds
under two additional conditions in the other case that $V$ is bounded.
In the special case that the sequence of continuous maps  $\{f_n\}_{n=0}^{\infty}$  converges to some continuous map $f$, some chaotic
 behavior of system (1.1) is investigated by using that of system (1.2) in Section 4.

\bigskip

\noindent{\bf 2. Preliminaries }\medskip

In this section, some basic concepts and lemmas are given. Especially, similarities and differences between  properties of the three conditions of Devaney chaos for system (1.1) and those for system (1.2) are presented.

For convenience, denote $f_{0,\infty}:=\{f_n\}_{n=0}^{\infty}$,  $f_n^{i}:=f_{n+i-1}\circ\cdots\circ f_n$, and $f_n^{0}:=id$. For a fixed $x_0\in D_0$, set $O(x_0):=\{f_0^i(x_0):i\geq 0\}$. Let $x\in X$ and $\varepsilon>0$.  By $B_{\varepsilon}(x)$ and $\overline{B}_{\varepsilon}(x)$  denote the open and closed  balls of radius $\varepsilon$, centered at $x$, respectively. \medskip

\noindent{\bf Definition 2.1} [21, Definition 2.1].  A point $p\in D_0$ is called $k$-periodic for system (1.1) if $f_0^{n+k}(p)=f_0^{n}(p), n\geq0$.
Further, $k$ is called the prime period of $p$ if $f_0^{i}(p) \neq p$, $1\leq i\leq k-1$.  In the special case of
$k=1$, i.e., $f_n(p)=p$ for all $n\geq0$, $p$ is called a fixed point of system (1.1). The set of all the periodic points
is denoted by $P(f_{0,\infty})$.
\medskip

\noindent{\bf Remark 2.1.} It is easy to see that two  periodic orbits of system (1.2) either do not intersect or agree completely.  However, it is not true for system (1.1) in general. See the following example: \medskip

\noindent{\bf Example 2.1.} Let $X=\{a, b, c, d\}$ be equipped with discrete metric. Set $D_{4n}=\{a,b\}$, $D_{4n+1}=\{b,c\}$, $D_{4n+2}=\{a,d\}$, $D_{4n+3}=\{a,b\}$, $f_{2n}(a)=b$, $f_{2n+1}(b)=a$, $f_{4n}(b)=c$,  $f_{4n+1}(c)=d$,  $f_{4n+2}(d)=a$,  and $f_{4n+3}(a)=b$ for each $n\geq0$.
Then $a$ is a $2$-periodic point and $O(a)=\{a,b\}$; $b$ is a $4$-periodic point and $O(b)=\{b, c, d, a\}$. It is evident that $O(a)\cap O(b)=\{a,b\}$ and $O(a)\neq O(b)$.\medskip

Let $V$ be a subset of $D_0$. For any two nonempty relatively open subsets $U_0$ and $V_0$ of  $V$, denote
$N(U_0,V_0):=\{n\in \mathbb{N_{+}}:f_0^{n}(U_0)\cap  V_0\neq\emptyset\}.$
\medskip

\noindent{\bf Definition 2.2} [21, Definition 2.2]. Let $V$ be a nonempty subset of $D_0$. System (1.1) is said to be topologically
transitive in $V$ if $N(U_0,V_0)\neq\emptyset$ for any two nonempty relatively open subsets $U_0$ and $V_0$ of $V$.\medskip

\noindent{\bf Remark 2.2.} For system (1.2), it is easy to show  that if  it is topologically
transitive in a finite set $V$, then $V$ consists of a periodic orbit.
However, it is not true for system (1.1) in general. \medskip

\noindent{\bf Example 2.2.}  Let $X=\{a, b\}$ be equipped with discrete metric. Set $V=D_n=X$ for each $n\geq 0$, $f_{0}(a)=f_{1}(a)=b$, $f_{0}(b)=f_{1}(b)=a$,  $f_{m}(a)=a$, and $f_{m}(b)=b$ for each $m\geq2$.
Then system (1.1) is  topologically
transitive in $V$. However, there are no periodic points in $V$. \medskip

It is known   that  if  system (1.2) is topologically
transitive in an infinite set $V$, then $V$ contains no isolated points  [13].
 For  system (1.1), a similar result is not true in general.\medskip

\noindent{\bf Example 2.3.} Let $X=[0,1]\cup\{2\}$ be equipped with the induced topology from the Euclidean space $\mathbb{R}$ and $\{a_n\}_{n=0}^{\infty}$ be the sequence of rational numbers in $[0,1]$. Set $V=D_n=X$, $f_{2n}(x)=a_n$, and $f_{2n+1}(x)=2$ for each $n\geq0$ and $x\in V$. Though system (1.1) is topologically
transitive in the infinite set $V$,  $V$ contains an isolated point $2$.\medskip

 However, we have the following result, which will be useful in the sequent sections:\medskip

\noindent{\bf Lemma 2.1.} {\it  Let $V$ be an infinite subset of  $D_0$.
If system {\rm(1.1)} is topologically transitive in $V$  and $P(f_{0,\infty})\cap V$ is dense in $V$,
then there are no isolated points in $V$. }\medskip

\noindent{\bf Proof.} Suppose that there exists an isolated point  $x_0\in V$. Then there exists $\varepsilon_0>0$ such that $B_{\varepsilon_0}(x_0)\cap V=\{x_0\}$. Since $P(f_{0,\infty})\cap V$ is dense in $V$, $x_0$ is a periodic point. So one can choose  $y_0\in V\backslash O(x_0)$
due to the fact that $V$ is infinite. Let $\delta_0:=d (y_0, O(x_0))>0$. Then  $f_0^{n}(B_{\varepsilon_0}(x_0)\cap V)\cap (B_{\delta_0/3}(y_0) \cap V)=\emptyset$ for each $n\geq 0$, which contradicts  the assumption that system (1.1) is topologically transitive in $V$.
This completes the proof.\medskip

\noindent{\bf Definition 2.3} [21, Definition 2.3]. Let $V$ be a nonempty subset of $D_0$. System (1.1) is said to have sensitive
dependence on initial conditions (briefly, system (1.1) is sensitive) in $V$ if there exists a constant $\delta_0> 0$  such that for any
$x_0\in V$ and any neighborhood $U$ of $x_0$, there exist $y_0\in V\cap U$ and a positive integer $n$
such that $d(f_0^n(x_0), f_0^n (y_0))>\delta_0$. The constant $\delta_0$ is called a sensitivity constant of system (1.1)
in $V$.\medskip

\noindent{\bf Remark 2.3.} It is easy to see that if there exists an isolated point in $V$, then system (1.1) is not sensitive in $V$.\medskip

\noindent{\bf Definition 2.4} [21, Definition 2.5]. Let $V$ be a nonempty subset of $D_0$. System (1.1) is said to be chaotic in the
sense of Devaney  on $V$ if\begin{itemize}\vspace{-0.2cm}
\item[{\rm (1)}] system (1.1) is topologically transitive in $V$;\vspace{-0.2cm}
\item[{\rm (2)}] $P(f_{0,\infty})\cap V$ is dense in $V$;\vspace{-0.2cm}
\item[{\rm (3)}] system (1.1) has sensitive dependence on initial conditions in $V$.\vspace{-0.15cm}
\end{itemize}

The next lemma shows that $(1)$ and $(2)$ imply $(3)$ of Devaney chaos  for  a continuous map.\medskip

\noindent{\bf Lemma 2.2} [3]. Let $X$ be a metric space and infinite, and $f:X\to X$ be a continuous map. If $f$ is topologically transitive and has dense periodic points in $X$, then
$f$ is sensitive in $X$.\medskip

\noindent{\bf Definition 2.5} [21, Definition 2.4]. If $V$  is a nonempty subset of $\bigcap_{n=0}^{\infty} D_n\neq\emptyset$ satisfying  $f_0^n(V)\subset V$ for all
$n\geq0$,  then $V$ is called a totally invariant set of system $(1.1)$. \medskip

\noindent{\bf Definition 2.6} [22, Definition 2.6].  Assume that $D_n$ and $E_n$ are two subsets of a metric space $(X, d)$,
and  $h_n: D_n \to E_n$
is a  map for each $n \geq 0$. The sequence of maps $\{h_n\}_{n=0}^{\infty}$
is said to
be equi-continuous in $\{D_n\}_{n=0}^{\infty}$ if for any $\varepsilon>0$, there exists a constant $\delta>0$ such that
$d(h_n(x), h_n(y))<\varepsilon$ for all $x, y \in D_n$ with $d(x, y) < \delta$ and for all $n \geq 0$.\medskip

\noindent{\bf Definition 2.7.} An increasing sequence of nonnegative integers $\{n_k\}_{k=1}^{\infty}$
is said to be syndetic if there exists an integer $l \geq 1 $ such that
$n_{k+1}-n_k \leq l$ for all $k\geq  1.$\medskip

  Now we introduce the concept of topological ergodicity for system (1.1).\medskip

\noindent{\bf Definition 2.8.}  Let $V$ be a nonempty subset of $D_0$. System (1.1)
is said to be topologically
ergodic in $V$ if $N(U_0,V_0)$ is syndetic for any two nonempty relatively open subsets $U_0$ and $V_0$ of $V$.\medskip

It was shown that    system (1.2) is topologically ergodic if it
satisfies (1) and (2) in Definition 1.1  (See [1] or [9]).
The next result shows that a similar result is true in the non-autonomous case.\medskip

\noindent{\bf Theorem 2.1.} {\it Let  $V$ be a totally invariant set of system $(1.1)$.
Assume that $f_n$ is  continuous in $V$ for each $n\geq 0$. If system {\rm(1.1)} is topologically transitive in $V$ and $P(f_{0,\infty})\cap V$ is dense in $V$,
then for any two nonempty relatively open subsets $U_0$ and $V_0$ of $V$,
there exists $p_0\in U_0\cap P(f_{0,\infty})$ such that $O(p_0)\cap V_0\neq\emptyset$.
Consequently, system {\rm(1.1)} is topologically ergodic in $V$. }\medskip

\noindent{\bf Proof.} Let $U_0$ and $V_0$ be two nonempty relatively open  subsets of $V$.  Since system (1.1)
is topologically transitive in $V$, there exist $x_0\in U_0$ and $n_0> 0$ such that $f_0^{n_0}(x_0)\in V_0$. Fix  $\varepsilon>0$ satisfying that $B_{\varepsilon}(f_0^{n_0}(x_0))\cap V\subset V_0$.
By the continuity of  $f_0^{n_0}$  in $V$ and the assumption that $V$ is a totally invariant set of system (1.1), there exists $\delta>0$ such that $B_{\delta}(x_0)\cap V \subset U_0$ and
$f_0^{n_0}(B_{\delta}(x_0)\cap V)\subset B_{\varepsilon}(f_0^{n_0}(x_0))\cap V\subset V_0$.
 Since $P(f_{0,\infty})\cap V$ is dense in $V$, there exists
$p_0\in B_{\delta}(x_0)\cap V\cap P(f_{0,\infty})\subset U_0$ with period $m_0$.  Then $f_0^{n_0+km_0}(p_0)=f_0^{n_0}(p_0)\in V_0$ for each $k\geq0$.
Thus, $O(p_0)\cap V_0\neq\emptyset$ and $\{ n_0+km_0\}_{k=0}^{\infty}\subset N(U_0,V_0)$,
which implies  that $N(U_0,V_0)$ is syndetic.
Hence, system (1.1) is topologically ergodic in $V$. This completes the proof.\bigskip

\noindent{\bf 3. Transitivity and density of periodic points imply sensitivity}\medskip

In this section, we shall  show that transitivity and density of periodic points imply
sensitivity  in $V$ for system (1.1)  under certain conditions. The discussions are divided
into the two cases:  $V$ is unbounded and bounded.
\medskip

In the case that $V$ is unbounded, we get the following result:\medskip

\noindent{\bf Theorem 3.1.} {\it Let $V$ be an unbounded subset of  $D_0$.
If system {\rm(1.1)} is topologically transitive in $V$ and $P(f_{0,\infty})\cap V$ is dense in $V$, then it is sensitive in $V$ and the sensitivity constant can be an arbitrarily positive number.}\medskip

\noindent{\bf Proof.} Let $\delta>0$ and $a\geq2$ be two arbitrarily given number, $x$ be an arbitrary point in $V$, and $U(x)$ be any  neighborhood of $x$ in $X$. Then we show that there exist $y_0\in U(x)\cap V$
and an integer $m_0>0$ such that
\vspace{-0.2cm}$$d(f_0^{m_0}(x),f_0^{m_0}(y_0))> (a-1)\delta/2.\vspace{-0.2cm}$$

Since system (1.1) is topologically transitive in $V$ and $P(f_{0,\infty})\cap V$
is dense in $V$, by Lemma 2.1
there are no isolated points in $V$ and thus there exists a periodic point
$q\in U(x)\cap V$ such that  $q\neq x$. Set
\vspace{-0.2cm}$$L_{x,q}:=\max_{z\in O(q)} {d(x,z)}>0.                                   \eqno(3.1)\vspace{-0.2cm}$$
Then the discussions are divided into two cases.

 $\rm{\bf{Case \;1.}}$ $L_{x,q}\geq\delta$.

 By (3.1), $O(q)\subset\overline {B}_{L_{x,q}}(x)$.
Since $V$ is unbounded, $\left(X\backslash\overline{B}_{ aL_{x,q}}(x)\right)\cap V$ is a nonempty open subset of $V$.
By the transitivity of system (1.1) in $V$, there exist $y_1\in U(x)\cap V$
and an integer $m_1>0$ such that $f_0^{m_1}(y_1)\in \left(X\backslash\overline{B}_{aL_{x,q}}(x)\right)\cap V$.
Then
\vspace{-0.2cm}$$d(f_0^{m_1}(q),f_0^{m_1}(y_1))\geq d(x,f_0^{m_1}(y_1))-d(x,f_0^{m_1}(q))> aL_{x,q}-L_{x,q}\geq(a-1)\delta.\vspace{-0.2cm}$$
Thus, one gets that
\vspace{-0.2cm}$${\rm either}\; \;d(f_0^{m_1}(x),f_0^{m_1}(y_1))> (a-1)\delta/2\;\; {\rm or}\;\;
d(f_0^{m_1}(x),f_0^{m_1}(q))> (a-1)\delta/2.\eqno(3.2)\vspace{-0.2cm}$$

 $\rm{\bf{Case \;2.}}$ $L_{x,q}<\delta$.

Since $V$ is unbounded, $\left(X\backslash\overline{B}_{ a\delta}(x)\right)\cap V$ is a nonempty open subset of $V$.
As system (1.1) is topologically transitive in $V$, there exist $y_2\in U(x)\cap V$
and an integer $m_2>0$ such that $f_0^{m_2}(y_2)\in \left(X\backslash\overline{B}_{a\delta}(x)\right)\cap V$.
This, together with the fact that $f_0^{m_2}(q)\in \overline {B}_{L_{x,q}}(x)\subset\overline {B}_{\delta}(x),$ implies that
\vspace{-0.2cm}$$d(f_0^{m_2}(q),f_0^{m_2}(y_2))\geq d(x,f_0^{m_2}(y_2))-d(x,f_0^{m_2}(q))> a\delta-\delta=(a-1)\delta.\vspace{-0.2cm}$$
Thus,
\vspace{-0.2cm}$${\rm either}\;\; d(f_0^{m_2}(x),f_0^{m_2}(y_2))> (a-1)\delta/2 \;\;{\rm or}\;\;
d(f_0^{m_2}(x),f_0^{m_2}(q))> (a-1)\delta/2. \eqno(3.3)\vspace{-0.2cm}$$

It follows from (3.2) and (3.3) that system (1.1) is sensitive in $V$ with sensitivity constant $(a-1)\delta/2$.
Since $a\geq 2$ and $\delta>0$ are arbitrary, the sensitivity constant of system (1.1) in $V$ can be any positive number.
The proof is  complete.\medskip

The following example is given to illustrate Theorem 3.1.\medskip

\noindent{\bf Example 3.1.} Let $I$ be an unbounded interval in $\mathbb{R}$ and $Q:=\{a_n\}_{n=0}^{\infty}$ denote the sequence of rational numbers in $I$.  For each $n\geq 0$, set  $V=D_n=X:=I$, $f_n(x)=a_n$ for any $x\in Q$, and $f_n(x)=x$ for any $x\in X\backslash Q$. Then $P(f_{0,\infty})=V\backslash Q$ is dense in $V$. Let $U_0$ and $V_0$ be any two nonempty open subsets of $V$ and fix any $x_0\in U_0\cap Q$. Since there exists $n_0\geq 0$ such that $a_{n_0}\in V_0$,  thus $f_0^{n_0+1}(x_0)=a_{n_0}\in V_0$ and hence system (1.1) is topologically transitive in $V$. By Theorem 3.1, system (1.1) is sensitive in $V$ and its sensitivity constant can be an arbitrarily positive number.  \medskip

Now, we consider the other case that $V$ is bounded. The following lemma is needed.\medskip

\noindent{\bf Lemma 3.1.} {\it Let $V$ be an infinite subset of $D_0$ and $A\subset V$ be dense in $V$. If system {\rm(1.1)} is not sensitive
in $V$, then for any given $\delta>0$, there exist  $y_0\in A$ and $\varepsilon_0\in(0,\delta)$ such that
\vspace{-0.2cm}$$f_0^n(B_{\varepsilon_0}(y_0)\cap V)\subset B_{\delta}(f_0^{n}(y_0)), \;n\geq 0.\eqno(3.4)\vspace{-0.2cm}$$
}
\noindent{\bf Proof.}  Fix any $\delta>0$. Since system {\rm(1.1)} is not sensitive
 in $V$,   there exist $x_0\in V$ and $\varepsilon_0\in(0,\delta)$ such that $d(f_0^n(y),f_0^n(x_0))<\delta/2$ for any  $y\in B_{2\varepsilon_0}(x_0)\cap V$ and each $n>0$. Since $A$ is dense in $V$, there exists  $y_0\in  A $ such that $y_0\in B_{\varepsilon_0}(x_0)$.  Then for any $y\in B_{\varepsilon_0}(y_0)\cap V$, $d(y,x_0)\leq d(y,y_0) +d(y_0,x_0)<2\varepsilon_0$ and thus $d(f_0^n(y), f_0^n(x_0))<\delta/2$ for each $n>0$. Hence, $d(f_0^n(y),f_0^n(y_0))\leq
d(f_0^n(y),f_0^n(x_0))+d(f_0^n(x_0),f_0^n(y_0))<\delta/2+\delta/2=\delta$ for any $y\in B_{\varepsilon_0}(y_0)\cap V$ and each $n>0$.
This, together with the fact that $0<\varepsilon_0<\delta$,  implies that
(3.4) holds.  The proof is complete.\medskip

\noindent{\bf Theorem 3.2.} {\it Let  $V$ be a bounded and infinite subset of $D_0$. Assume that $\{f_n\}_{n=0}^{\infty}$ is equi-continuous in $\{D_n\}_{n=0}^{\infty}$ and  system {\rm(1.1)} has a fixed point in $V$.
If system {\rm(1.1)} is topologically transitive in $V$ and $P(f_{0,\infty})\cap V$ is dense in $V$, then it is sensitive in $V$.}\medskip

\noindent{\bf Proof.} Set $\delta=d(V)/7$, where $d(V):=\sup_{x,y\in V}d(x,y)>0$. Suppose that system {\rm(1.1)} is not sensitive in $V$. Then by Lemma 3.1 and the assumption that $P(f_{0,\infty})\cap V$ is dense in $V$, there exist $p_0\in P(f_{0,\infty})\cap V$ with prime period $n_0$ and $0<\varepsilon_0<\delta$ such that
\vspace{-0.2cm}$$f_0^n(B_{\varepsilon_0}(p_0)\cap V)\subset B_{\delta}(f_0^{n}(p_0)), \;n\geq 0.\eqno(3.5)\vspace{-0.2cm}$$
Note that there are no isolated points in $V$ by Lemma 2.1.   Since system {\rm(1.1)} is topologically transitive in $V$ and $p_0$ is a $n_0$-periodic point in $V$, by $(3.5)$ one gets that
\vspace{-0.2cm}$$V\subset\overline{\bigcup_{n=0}^{\infty}f_0^n(B_{\varepsilon_0}(p_0)\cap V)}\subset\overline{\bigcup_{n=0}^{\infty}B_{\delta}(f_0^n(p_0))}=\overline{\bigcup_{n=k_0}^{k_0+n_0-1}B_{\delta}(f_0^n(p_0))}\eqno(3.6)\vspace{-0.2cm}$$
for any fixed $k_0\geq0$.
Since $\{f_n\}_{n=0}^{\infty}$ is equi-continuous in $\{D_n\}_{n=0}^{\infty}$,
it can be easily verified that there exists $0<\varepsilon_1<\delta$  such that for
any $n \geq 0$ and any $x, y\in D_n$ with $d(x,y)<\varepsilon_1$,
\vspace{-0.2cm}$$d(f_n^i(x),f_n^i(y))<\delta,\; 0\leq i\leq n_0-1.\eqno(3.7)\vspace{-0.2cm}$$
Let $q_0\in V$  be the fixed point of system (1.1). Then $f_n^i(q_0)=q_0$, and thus by (3.7),
\vspace{-0.2cm}$$f_n^i(B_{\varepsilon_1}(q_0)\cap D_n)\subset B_{\delta}(q_0), \;n\geq 0,\;0\leq i\leq n_0-1.\eqno(3.8)\vspace{-0.2cm}$$
 Since system (1.1) is topologically transitive in $V$,  there exists an integer $m_0>0$ such that
$f_0^{m_0}(B_{\varepsilon_0}(p_0)\cap V)\cap (B_{\varepsilon_1}(q_0)\cap V)\neq\emptyset. $
 Hence, there exists $y_0\in B_{\varepsilon_0}(p_0)\cap V$ such that
\vspace{-0.2cm}$$f_0^{m_0}(y_0)\in B_{\varepsilon_1}(q_0)\cap V\cap D_{m_0}.\eqno(3.9)\vspace{-0.2cm}$$
By (3.8) and (3.9), one gets that
\vspace{-0.2cm}$$f_0^{m_0+i}(y_0)=f_{m_0}^i\circ f_0^{m_0}(y_0)\in f_{m_0}^i(B_{\varepsilon_1}(q_0)\cap D_{m_0})\subset B_{\delta}(q_0),\;0\leq i \leq n_0-1.\eqno(3.10)\vspace{-0.2cm}$$
 Since $y_0\in B_{\varepsilon_0}(p_0)\cap V$, by (3.5) one has that
 \vspace{-0.2cm}$$f_0^{m_0+i}(y_0)\in f_0^{m_0+i}(B_{\varepsilon_0}(p_0)\cap V)\subset B_{\delta}(f_0^{m_0+i}(p_0)),\;0\leq i \leq n_0-1.\eqno(3.11) \vspace{-0.2cm}$$
 By (3.10) and (3.11), one gets that for any given $0\leq i\leq n_0-1$  and any $z_i\in B_{\delta}(f_0^{m_0+i}(p_0))$,
  \vspace{-0.4cm}$$d(z_i,q_0)\leq d(z_i,f_0^{m_0+i}(p_0))+d(f_0^{m_0+i}(p_0),f_0^{m_0+i}(y_0))+d(f_0^{m_0+i}(y_0),q_0)<3\delta.\vspace{-0.2cm}$$
  This implies that
 \vspace{-0.2cm}$$\bigcup_{i=0}^{n_0-1}B_{\delta}(f_0^{m_0+i}(p_0))\subset B_{3\delta}(q_0),\vspace{-0.2cm}$$
 and thus by (3.6) one gets that
   \vspace{-0.2cm}$$V\subset\overline{\bigcup_{n=m_0}^{m_0+n_0-1}B_{\delta}(f_0^{n}(p_0))}\subset \overline{B}_{3\delta}(q_0).\vspace{-0.2cm}$$
   Hence, $d(V)\leq6\delta$, which  contradicts the fact that $d(V)=7\delta$. The proof is complete.
   \medskip

   The following example is given to illustrate Theorem 3.2.\medskip

\noindent{\bf Example 3.2.} Let $X=D_{n}=[0,1]$,   $f_{2n}(x)=x$,  and
\vspace{-0.1cm}$$\vspace{-0.2cm}
f_{2n+1}(x) =\begin{cases} 2x & \text{ if } x\in [0,1/2], \\
                     2(1-x) & \text{ if } x\in(1/2,1], \end{cases}
\vspace{-0.05cm}$$
be the tent map for each $n\geq0$.  Set $V=D_0$. Since $f_{2n+1}\circ f_{2n}$
is also the tent map in $[0,1]$, it is easy to check that system (1.1)
satisfies all the conditions in Theorem 3.2. Thus, it is sensitive in $V$.   \medskip

For any given positive integer $N$, denote
    \vspace{-0.2cm}$$P_N(f_{0,\infty}):=\{p: p\in P(f_{0,\infty}) \; {\rm with\; prime\; period \;no\; more\; than}\; N\}. \vspace{-0.2cm}$$
It is evident that the condition that $P_N(f_{0,\infty})\cap V$ is dense in $V$ for some $N>0$
is stronger than condition $(2)$ in Definition 2.4.\medskip

\noindent{\bf Theorem 3.3.} {\it  Let $V$ be a bounded and infinite subset of  $D_0$.
If system {\rm(1.1)} is topologically transitive in $V$ and $P_N(f_{0,\infty})\cap V$ is dense in $V$ for some $N>0$, then it is sensitive
 in $V$.}\medskip

\noindent{\bf Proof.} Since $V$ is infinite, there exist  $N+1$ different points  $x_i\in V, 0\leq i\leq N$, and $\delta>0$ such that  $B_{3\delta}(x_i), 0\leq i\leq N$, are mutually disjointed.  Suppose that system {\rm(1.1)} is not sensitive in $V$. Then with a similar argument to the proof of  (3.6), one  gets that there exists $p_1\in P_N(f_{0,\infty})\cap V$ with prime period $n_1\leq N$ such that
\vspace{-0.2cm}$$V\subset\bigcup_{n=0}^{n_1-1}\overline{B}_{\delta}(f_0^n(p_1)).\eqno(3.12)\vspace{-0.2cm}$$

By noting that $B_{3\delta}(x_i), 0\leq i\leq N$, are mutually disjointed, it is easy to verify that for each  $0\leq j\leq n_1-1$, $\overline{B}_{\delta}(f_0^j(p_1))$ contains at most one element of $\{x_i\}_{i=0}^{N}$.
Thus there exists some $i_0$ such that $0\leq i_0\leq N$ and $x_{i_0}\notin \cup_{n=0}^{n_1-1}\overline{B}_{\delta}(f_0^n(p_1))$. This contradicts $(3.12)$. Hence, system (1.1) is sensitive in $V$. The proof is complete.\medskip

\noindent{\bf Example 3.3.} Replace the unbounded interval $I$ in $\mathbb{R}$ by a bounded one in Example 3.1. Then  system (1.1) in Example 3.1 is sensitive in $V$ by Theorem 3.3. \bigskip

\noindent{\bf 4. The special case that $\{f_n\}_{n=0}^{\infty}$  converges to $f$}\medskip

In this section, we shall first study relationships between some chaotic behavior of system (1.1) and  that of system (1.2) if the sequence of continuous maps $\{f_n\}_{n=0}^{\infty}$  converges to a continuous map $f$ and condition $(2)$  in Definition 2.4 holds, and  then show that transitivity and density of periodic points imply sensitivity for system (1.1) in this special case. In addition, if $f$ (or $f_0$) is surjective, then $f_n=f$ for each $n\geq 0$.

We first consider the relationship between periodic points of system (1.1) and those of system (1.2).\medskip

\noindent{\bf Lemma 4.1.} {\it Let  $V$ be  a totally invariant set of system $(1.1)$ and
 $\{f_n\}_{n=0}^{\infty}$ converge to $f$  in $V$. If $x_0\in V$ is a $k_0$-periodic point of system {\rm(1.1)}, then
$f_{0}^{i}(x_0)=f^{i}(x_0)$ for each  $i\geq0$.  Consequently, $x_0$ is a $k_0$-periodic point of system $(1.2)$.}\medskip

\noindent{\bf Proof.}
Let $x_0\in V$ be a $k_0$-periodic point of system (1.1) and $i\geq0$.  Then \vspace{-0.2cm}$$f_0^{i+1}(x_0)=f_0^{nk_0+i+1}(x_0)=f_{nk_0+i}(f_0^{nk_0+i}(x_0))=f_{nk_0+i}(f_0^i(x_0)), \;n\geq0.\vspace{-0.2cm}$$
 Since $V$ is  a totally invariant set of system $(1.1)$ and $f_{nk_0+i}\to f$  in $V$ as $n\to\infty$,  $f_0^{i+1}(x_0)=f(f_0^i(x_0))$. By induction, one gets that $f_{0}^{i}(x_0)=f^{i}(x_0)$ for each $i\geq0$. This completes the  proof.\medskip

Now, we study  the relationship between some chaotic  behavior of system (1.1) and that of system (1.2).\medskip

\noindent{\bf Proposition 4.1.} {\it Let $V$ be  a totally invariant set of system $(1.1)$
and $\{f_n\}_{n=0}^{\infty}$ converge to a map $f$  in $V$. If $f_n$ for all $n\ge 0$ and
$f$ are continuous in $V$, and  $P(f_{0,\infty})\cap V$ is dense in $V$, then
 \begin{itemize}\vspace{-0.2cm}
\item[{\rm (i)}] $f_0^n=f^n$ in $V$ for each $n\geq1$;\vspace{-0.2cm}
\item[{\rm (ii)}]  topological transitivity of system $(1.1)$ in $V$  is equivalent to that of system $(1.2)$ in $V$;\vspace{-0.2cm}
\item[{\rm (iii)}] sensitivity of system $(1.1)$ in $V$ is equivalent to that of system $(1.2)$ in $V$.\vspace{-0.2cm}
 \end{itemize}}

\noindent{\bf Proof.}  First, we show that (i) holds. Fix any $x\in V$. Since $P(f_{0,\infty})\cap V$ is dense in $V$, there exists a sequence $\{p_m\}_{m=0}^{\infty}\subset P(f_{0,\infty})\cap V$ such that $p_m\to x$ as $m\to\infty$.  By Lemma 4.1, we have that
\vspace{-0.2cm}$$f_0^n(p_m)=f^n(p_m),\; m\geq 0,\;n\geq 1. \eqno(4.1)\vspace{-0.2cm}$$
 Letting $m\to\infty$ in (4.1) for any given $n\geq 1$, one has that
 $f_0^n(x)=f^n(x)$ by the continuity of $f_0^n$ and $f^n$. Thus assertion  (i) holds.

 Assertions (ii) and (iii) can be shown by (i) immediately. This completes the proof.\medskip

Next, we shall show that  the two former conditions imply the last one in Definition 2.4 by Lemma 2.2
and Proposition 4.1 under some assumptions.\medskip

\noindent{\bf Theorem 4.1.} {\it Let $V$ be  a totally invariant set of system $(1.1)$ and infinite.
Assume that $\{f_n\}_{n=0}^{\infty}$ converges to a map $f$  in $V$, and $f_n$ for all $n\ge 0$ and
$f$ are continuous in $V$.
If  system $(1.1)$ is topologically transitive in $V$ and $P(f_{0,\infty})\cap V$ is dense in $V$, then it is sensitive  in $V$.}\medskip

\noindent{\bf Proof.}  Since system $(1.1)$ is topologically transitive in $V$, system (1.2) is also topologically transitive in $V$ by (ii) of Proposition 4.1. Because  $P(f_{0,\infty})\cap V$ is dense in $V$,  it follows from Lemma 4.1 that $P(f)\cap V$ is dense in $V$. By Lemma 2.2, system (1.2) is sensitive in $V$.
Therefore, system (1.1) is sensitive in $V$ by (iii) of Proposition 4.1. The proof is complete. \medskip

\noindent{\bf Remark 4.1.} By Lemma 4.1 and Proposition 4.1, Devaney chaos of system (1.1) implies that of system (1.2) under  the assumptions of Theorem 4.1. However, the converse is not true. See the following example:\medskip

\noindent{\bf Example 4.1.} Let $V=D_n=X:=[0,1]$,
\vspace{-0.2cm}$$\begin{array}{cccc}\begin{array}{cccc}
f_0(x) =\begin{cases} 2x & \text{ if } x\in [0,1/3], \\
                     2/3 & \text{ if } x\in(1/3,2/3), \\
                      x & \text{ if } x\in [2/3,1],\end{cases} \end{array}
&\begin{array}{cccc}
f=f_n(x) =\begin{cases} 2x & \text{ if } x\in [0,1/2], \\
                   -2x+2 & \text{ if }  x\in (1/2,1], \end{cases}
\end{array}\end{array}\vspace{-0.2cm}$$
where $n\geq 1$. Then $\{f_n\}_{n=0}^{\infty}$  converges  to  $f$ in $V$. Let $U_0=(1/3,2/3)$ and $V_0=(0,1/3)$. It is clear that $f_0^n(U_0)=\{2/3\}$ and hence $f_0^n(U_0)\cap V_0=\emptyset$ for each $n> 0$. This implies that system (1.1) is not chaotic in the sense of Devaney.  However, $f$ is the tent map in $[0,1]$ and hence system (1.2) is chaotic in the sense of Devaney.\medskip

The following result shows that $f_n=f$  in $V$ for each $n\ge 0$ if all the conditions in Proposition 4.1 are satisfied and in addition $f_0$  is surjective in $V$. \medskip

 \noindent{\bf Proposition 4.2.} {\it Let $V$ be  a totally invariant set of system $(1.1)$
 and $f_n$ be continuous in $V$ for each $n\geq0$.  Assume that $P(f_{0,\infty})\cap V$
 is dense in $V$ and $f_0$ is surjective in $V$.
Then  $\{f_n\}_{n=0}^{\infty}$ converges to a  continuous  map $f$  in $V$
if and only if $f_n=f$ in $V$ for each $n\geq 0$.}\medskip

 \noindent{\bf Proof.} It suffices to show the necessity.
  By (i) of Proposition 4.1, $f_0=f$ in $V$. Thus $f$ is surjective in $V$.
Suppose that for $k\geq 0$, $f_j=f$ in $V$ for each $0\leq j\leq k$. Then $f_j,0\leq j\leq k$, are surjective in $V$. For any $y\in V$, there exist $x_j\in V$, $0\leq j \leq k$, such that $y=f_k(x_k)=f(x_k)$, $x_{i+1}=f_i(x_i)=f(x_i)$, $0\leq i\leq k-1$. It follows from (i) of Proposition 4.1  that $f_{k+1}(y)=f_0^{k+2}(x_0)=f^{k+2}(x_0)=f(y)$, which implies that $f_{k+1}=f$ in $V$. Hence, the necessity holds by induction.
This completes the proof.
\medskip

The following result is a direct consequence of  Propositions 4.1 and 4.2.\medskip

\noindent{\bf Corollary 4.1.} {\it Let all the conditions  in Proposition {\rm 4.1} hold.  Assume that $f$ is surjective in $V$. Then $f_n=f$ for  each $n\geq 0$.}\medskip

 In the case that the metric space is compact, we have the following result: \medskip

 \noindent{\bf Proposition 4.3.} {\it Let $(X, d)$ be a compact  metric space and $f_n:X\to X$ be continuous
for each $n\geq0$. Assume that system $(1.1)$  is topologically transitive and  $P(f_{0,\infty})$ is dense in $X$.
Then the following statements are equivalent:
 \begin{itemize}\vspace{-0.2cm}
\item[{\rm (i)}] $f_n=f$ for each $n\in \mathbb{N}$;\vspace{-0.2cm}
\item[{\rm (ii)}] $\{f_n\}_{n=0}^{\infty}$ uniformly converges to a map $f$ in $X$;\vspace{-0.2cm}
\item[{\rm (iii)}] $\{f_n\}_{n=0}^{\infty}$ converges to a continuous map $f$ in $X$.\vspace{-0.2cm}
 \end{itemize} }

  \noindent{\bf Proof.} ${\rm(i)}\Rightarrow{\rm(ii)}$ and ${\rm(ii)}\Rightarrow{\rm(iii)}$ are obvious.

  ${\rm(iii)}\Rightarrow{\rm(i)}$: By Corollary 4.1, it suffices to show that $f$ is surjective in $X$.
  Otherwise, $f(X)\subsetneqq X$.  Since $f$ is continuous in $X$ and $X$ is compact, $f(X)$ is compact and thus  closed  in $X$. Hence, $X\backslash f(X)$ is a nonempty open subset in $X$. By (i) of Proposition 4.1, $f_0^n(U_0)=f^n(U_0)\subset f(X)$ for any fixed nonempty open subset $U_0$ and each $n> 0$. Thus, $f_0^n(U_0)\cap(X\backslash f(X))=\emptyset$ for each $n> 0$. This  contradicts  the assumption that system $(1.1)$ is topologically transitive in $X$. Therefore, $f$ is surjective in $X$. The proof is complete.
\end{document}